\documentclass[a4paper,11pt,reqno]{amsart}
\usepackage{amsmath, amsthm, amssymb, enumerate, mathtools, amsaddr}
\usepackage[a4paper, top=1in, bottom=1.25in, left=1.25in, right=1.25in]{geometry}
\usepackage{bbm}
\usepackage{enumitem}
\usepackage{hyperref}
\usepackage{amsrefs}

\newenvironment{proofclaim}[1][Proof of Claim]{\begin{proof}[#1]}{\end{proof}}

\def\eps{{\varepsilon}}

\begin{document}

\numberwithin{equation}{section}

\title[Upper density of monochromatic paths in edge-coloured $K_{\mathbb{N}}$ and $K_{\mathbb{N},\mathbb{N}}$]{Upper density of monochromatic paths in edge-coloured infinite complete graphs and bipartite graphs}
\date{\today}
\author{A. Nicholas Day}
\address[A. Nicholas Day]{Ume\r{a} University, Ume\r{a}, Sweden}
\email{a.nick.day@gmail.com}
\author{Allan Lo}
\address[Allan Lo]{University of Birmingham, Birmingham B15 2TT, United Kingdom}
\email{s.a.lo@bham.ac.uk}
\thanks{The research leading to these results was partially supported by EPSRC, grant no. EP/P002420/1 (A.~Lo).}

\newtheorem{Thm}{Theorem}[section]
\newtheorem{Lemma}[Thm]{Lemma}
\newtheorem{Def}[Thm]{Definition}
\newtheorem{Coro}[Thm]{Corollary}
\newtheorem{Conj}[Thm]{Conjecture}
\newtheorem{Prop}[Thm]{Proposition}
\newtheorem{Quest}[Thm]{Question}
\newtheorem{Claim}[Thm]{Claim}
\newtheorem{Case}{Case}
\newtheorem{Subcase}{Subcase}
\newtheorem{Case1}{Case}

\begin{abstract}
The upper density of an infinite graph~$G$ with $V(G) \subseteq \mathbb{N}$ is defined as $
\overline{d}(G) = \limsup_{n \rightarrow \infty}{|V(G) \cap \{1,\ldots,n\}|}/{n}$.
Let $K_{\mathbb{N}}$ be the infinite complete graph with vertex set~$\mathbb{N}$. 
Corsten, DeBiasio, Lamaison and Lang showed that in every $2$-edge-colouring of~$K_{\mathbb{N}}$, there exists a monochromatic path with upper density at least $(12 + \sqrt{8})/17$, which is best possible.
In this paper, we extend this result to $k$-edge-colouring of~$K_{\mathbb{N}}$ for $k \ge 3$. 
We conjecture that every $k$-edge-coloured~$K_{\mathbb{N}}$ contains a monochromatic path with upper density at least~$1/(k-1)$, which is best possible (when $k-1$ is a prime power).
We prove that this is true when $k = 3$ and asymptotically when $k =4$. 
Furthermore, we show that this problem can be deduced from its bipartite variant, which is of independent interest. 
\end{abstract}
\maketitle 

\section{Introduction}
Throughout the paper, a $k$-edge-colouring of a graph uses colours $1,2, \dots, k$. 
Given a $k$-edge-coloured~$K_n$, how long is the longest monochromatic path? 
This question is equivalent to asking for the $k$-colour Ramsey number of paths~$P_n$, denoted by~$R_k(P_n)$.
When $k=2$, Gerencs\'er and Gy\'arf\'as~\cite{GG} show that $R_2(P_n) = \lfloor 3n/2 \rfloor -1$.
When $k = 3$, Gy\'arf\'as, Ruszink\'o, S\'ark\"ozy and Szemer\'edi~\cite{GRSS} show that $R_3(P_n) = n +2 \lceil n/2 \rceil- 2$.
For $k \ge 4$, we know that $(k-1+o(1))n \le  R_k(P_n) \le (k-1/2+o(1))n$ by  Sun, Yang, Xu and Li~\cite{ramseypathlower} and Knierim and Su~\cite{ramseypathupper}, respectively. 
This implies that every $k$-edge-coloured~$K_n$ contains a monochromatic path of density between $2/(2k-1)$ and $1/(k-1)$. 
We ask the analogous question for infinite complete graphs. 

Let $K_{\mathbb{N}}$ be the infinite complete graph with vertex set~$\mathbb{N}$, where $\mathbb{N}$ is the set of strictly positive integers (i.e. without zero).
Given a set $A \subseteq \mathbb{N}$, the \textit{upper density of~$A$} is defined as 
\begin{equation}
\overline{d}(A) = \limsup_{n \rightarrow \infty}\frac{|A \cap [n]|}{n}, \nonumber
\end{equation}
where $[n] = \{1, \ldots, n\}$. 
Similarly, if $G$ is a graph with $V(G) \subseteq \mathbb{N}$, then the \textit{upper density of~$G$} is defined as $\overline{d}(G) = \overline{d}(V(G))$.  
Hence, the main focus of this paper is determine the largest upper density of a monochromatic path guaranteed in any $k$-edge-colouring of~$K_{\mathbb{N}}$.

A result of Rado~\cite{Rado} states that any $k$-edge-coloured~$K_{\mathbb{N}}$ can be partitioned into at most $k$ monochromatic paths. 
This implies that one of these paths must have upper density at least~$1/k$.
When $k=2$, Erd\H{o}s and Galvin~\cite{ErdosGalvin} proved that any $2$-edge-coloured~$K_{\mathbb{N}}$ contains a monochromatic path with upper density between $2/3$ and $8/9$.
The lower bound was further improved in~\cite{DeBiasioMchKenney, LoSMWang}.
Finally, Corsten, DeBiasio, Lamaison and~Lang~\cite{CDLL} proved that any $2$-edge-coloured $K_{\mathbb{N}}$ contains a monochromatic path with upper density at least $(12+\sqrt{8})/17$, which is best possible. 
See~\cite{BL, CDM, Lamaison} for other upper densities of monochromatic subgraphs in $2$-edge-coloured $K_{\mathbb{N}}$.

In this paper, we consider $k$-edge-colourings of $K_{\mathbb{N}}$ for $k \ge 3$; we conjecture that the picture is quite different from the case $k = 2$.
\begin{Conj}\label{conjecture:main-conjecture}
Let $k \ge 3$.
In any $k$-edge-colouring of~$K_{\mathbb{N}}$, there exists a monochromatic path~$P$ with $\overline{d}(P) \ge 1/({k-1})$.
\end{Conj}

We prove that this conjecture holds when $k = 3$ and asymptotically when $k =4$. 

\begin{Thm}\label{theorem:k-colourings-complete-graph}
In any $3$-edge-colouring of~$K_{\mathbb{N}}$, there exists a monochromatic path~$P$ with $\overline{d}(P) \ge 1/2$.
In any $4$-edge-colouring of~$K_{\mathbb{N}}$, $\sup \{ \overline{d}(P) : P \text{ is a monochromatic path} \} \ge 1/3$.
\end{Thm}

The following result shows that if Conjecture~\ref{conjecture:main-conjecture} holds, then it is sharp when $k-1$ is a prime power. 
In particular, Theorem~\ref{theorem:k-colourings-complete-graph} is best possible. 

\begin{Prop}[{\cite[Corollary~3.5]{DeBiasioMchKenney}}]\label{theorem:k-colourings-construction}
Let $k \ge 3$ and let $q$ be a prime power with $q \le k-1$. 
Then there exists a $k$-edge-colouring of $K_{\mathbb{N}}$ in which every monochromatic path~$P$ satisfies $\overline{d}(P) \le 1/q$.
\end{Prop}

\subsection{Monochromatic paths in complete bipartite infinite graphs.}

We now look at the bipartite variant of the above problem, which turns out to be closely related. 
Let $K_{\mathbb{N},\mathbb{N}}$ denote the set of all complete bipartite graphs on $\mathbb{N}$ where both vertex classes are infinite.
We will write $K_{V,W} \in K_{\mathbb{N},\mathbb{N}}$ to denote that $K_{V,W}$ is a graph in $K_{\mathbb{N},\mathbb{N}}$ such that $V$ and $W$ are both infinite, disjoint and together partition $\mathbb{N}$.

We investigate the largest guaranteed upper density of monochromatic paths in a $k$-edge-coloured $K_{V,W} \in K_{\mathbb{N},\mathbb{N}}$.
A result of D.~Soukup~\cite[Theorem 2.4.1]{Soukup} states that any $k$-edge-coloured~$K_{V,W}$ can be partitioned into at most $2k-1$ monochromatic paths. 
So one of them will have upper density at least $1/(2k-1)$.
By considering a $k$-edge-colouring of~$K_{V,W}$ derived from a proper $k$-edge-colouring of~$K_{k,k}$ (see Section~\ref{section:bipartite-work} for details), it is not difficult to see that $1/k$ is an upper bound on the upper density of a monochromatic path in~$K_{V,W}$.

\begin{Prop}\label{prop:bipartite-k-colourings-construction}
Let $K_{V,W} \in K_{\mathbb{N},\mathbb{N}}$.
For all $k \in \mathbb{N}$, there exists a $k$-edge-colouring of~$K_{V, W}$ in which every monochromatic path~$P$ satisfies $\overline{d}(P) \leqslant 1/k$.
\end{Prop}

We conjecture that this bound is in fact tight. 
\begin{Conj}\label{conjecture:bipartite-conjecture}
In any $k$-edge-colouring of $K_{V,W} \in K_{\mathbb{N},\mathbb{N}}$, there exists a monochromatic path~$P$ with $\overline{d}(P) \ge 1/k$.
\end{Conj}

The conjecture is trivial for $k=1$. 
In this paper, we show that the conjecture is true for $k =2$ and a weaker result when $k \ge 3$, which improves the lower bound of $1/(2k-1)$ coming from Soukup's result~\cite{Soukup}.

\begin{Thm}\label{theorem:multicolour-bipartite}
Let $K_{V,W} \in K_{\mathbb{N},\mathbb{N}}$. 
In any $2$-edge-colouring of~$K_{V,W}$, there exists a monochromatic path~$P$ with $\overline{d}(P) \ge 1/2$. 
In any $k$-edge-colouring of~$K_{V,W}$ with $k \ge 3$, $\sup \{ \overline{d}(P) : P \text{ is a monochromatic path} \} \ge 1/(2k-3)$.
\end{Thm}

The case $k=2$ for both Theorem~\ref{theorem:multicolour-bipartite} and Conjecture~\ref{conjecture:bipartite-conjecture} are also implied by a result of Corsten, DeBiasio and McKenney~\cite[Theorem~1.15]{CDM}.
They also conjectured a stronger version of Conjecture~\ref{conjecture:bipartite-conjecture}.

We show that one can deduce Conjecture~\ref{conjecture:main-conjecture} from Conjecture~\ref{conjecture:bipartite-conjecture} using the following result.

\begin{Thm}\label{theorem:main-theorem}
Let $k \ge 3$.
Let $\phi_k$ be such that, for all $K_{V,W} \in K_{\mathbb{N},\mathbb{N}}$ and all $k$-edge-colourings of $K_{V,W}$, there exists a monochromatic path~$P$ with $\overline{d}(P) \ge {\phi}_k$. 
 Then, in any $k$-edge-colouring of $K_{\mathbb{N}}$, there exists a monochromatic path~$P$ with $\overline{d}(P) \ge {\phi}_{k-1}$.
\end{Thm}

Therefore, Theorem~\ref{theorem:k-colourings-complete-graph} is a corollary of Theorems~\ref{theorem:multicolour-bipartite} and~\ref{theorem:main-theorem}.

A natural open problem is to show that Conjecture~\ref{conjecture:bipartite-conjecture} for all~$k \ge 3$. 
Note that when $k=3$, Theorem~\ref{theorem:main-theorem} shows that the supremum of the upper densities of monochromatic paths is at least~$1/3$.
It would be good to replace the supremum with the maximum.

\subsection{Notations and layout.}
Given a graph $G$, a vertex $v \in V(G)$, and a colouring of the edges of $G$, we write $N_{i}(v)$ for the set of vertices $x$ in $G$ such that edge $vx$ has colour~$i$.
A path of colour~$i$ refers to a monochromatic path of colour~$i$. 

The layout of this paper is as follows.  
In Section \ref{section:proof-of-main-theorem}, we employ some ideas from graph regularity to prove Theorem \ref{theorem:main-theorem}. 
In Section \ref{section:bipartite-work}, we investigate $k$-edge-coloured~$K_{\mathbb{N}, \mathbb{N}}$. 
In particular, we prove Proposition~\ref{prop:bipartite-k-colourings-construction} and Theorem~\ref{theorem:multicolour-bipartite}.

\section{Proof of Theorem \ref{theorem:main-theorem}}\label{section:proof-of-main-theorem}

We need a Ramsey-type result on matchings. 
Cockayne and Lorimer~\cite{CockayneLorimer} showed that in any $k$-edge-colouring of~$K_n$, there exists a monochromatic matching of size at least $\left\lfloor \frac{n+k-1}{k+1} \right\rfloor$.
We need the analogous result by~Omidi, Raesi and Rahimi~\cite{ORR2018}, which replaces $K_n$ with graphs with large minimum degree.

\begin{Thm}[{Omidi, Raesi and Rahimi~\cite[Corollary~1.4]{ORR2018}}]\label{lemma:regularity-ramsey}
Let $G$ be a graph on $n$ vertices with $\delta(G) \geq \frac{kn}{k+1}$.
Then, in any $k$-edge-colouring of~$G$, there exists a monochromatic matching of size at least $\left\lfloor \frac{n+k-1}{k+1} \right\rfloor$.
\end{Thm}

%

By a standard application of the regularity lemma, we can deduce that any large $k$-edge-coloured~$K_n$ contains a constant number of monochromatic paths covering about $2/(k+1)$ fraction of the vertices.
(Namely, we first apply the Szemer\'edi's regularity lemma to~$K_n$, then obtain a large monochromatic matching in the reduced graph using Theorem~\ref{lemma:regularity-ramsey}, and finally, convert the edges of the monochromatic matchings into long monochromatic paths by  the blow-up lemma.
For example, see the proof of~Lemma~3 in~\cite{FigajLuczak}.)
 

\begin{Coro} \label{cor:regularity-ramsey}
For all $k \ge 3 $ and all $\eps > 0$, there exists $ t_0 = t_0(k, \eps) > 0$ and $n_{0} = n_0(k, \eps) \in \mathbb{N}$ such that the following holds.  
In any $k$-edge-colouring of~$K_n$ with $n \ge n_0$, there exist vertex-disjoint paths $P_1, \dots, P_t$ such that $t \le t_0$, $\bigcup_{i \in [t]}P_i$ is monochromatic and $| \bigcup_{i \in [t]}P_i | \ge  \frac{(2-\eps)n}{k+1}$.
\end{Coro}

We now prove Theorem~\ref{theorem:main-theorem}.

\begin{proof}[Proof of Theorem~\ref{theorem:main-theorem}]
Fix a $k$-edge-colouring of~$K_{\mathbb{N}}$.  
Suppose to the contrary that there is no monochromatic path~$P$ with $\overline{d}(P) \ge \phi_{k-1}$. 
We say that a vertex $v \in V(K_{\mathbb{N}})$ is \emph{of colour~$i$} if $N_i(v)$ is infinite.  
Note that a vertex can have more than one colour. 

\begin{Claim} \label{clm:1}
For every finite set~$S \subseteq \mathbb{N}$, every $i \in [k]$ and every pair of vertices $x, y \in \mathbb{N} \setminus S$ of colour~$i$, there exists a path of colour~$i$ from $x$ to~$y$ that avoids~$S$. 
\end{Claim}

\begin{proofclaim}
Suppose  to the contrary that there exists a finite set $S \subseteq \mathbb{N}$, $i \in [k]$ and $x,y \in \mathbb{N} \setminus S$ of colour~$i$ such that there is no path from~$x$ to~$y$ of colour~$i$ which avoids~$S$. 
Let $X$ be the set of vertices that can be reached from $x$ by a path of colour $i$ that avoids~$S$, and let $Y = \mathbb{N} \setminus (X \cup S)$.
Since $N_i(x) \setminus S \subseteq X$ and $N_i(y) \setminus S \subseteq Y$,  each of $X$ and $Y$ is infinite. 
Moreover, there is no edge of colour~$i$ between $X$ and~$Y$.
Thus the infinite complete bipartite graph $K[X,Y]$ (with vertex classes $X$ and $Y$) is $(k-1)$-coloured.
So there exists a monochromatic path~$P$ in~$K[X,Y] \subseteq K_{\mathbb{N}}$ with 
\begin{align*}
\overline{d}(P) = \limsup_{n \rightarrow \infty}\frac{|V(P) \cap [n]|}{n} = \limsup_{n \rightarrow \infty}\frac{| (V(P) \cap [n]) \setminus S |}{|[n] \setminus S|} \ge \phi_{k-1},
\end{align*}
a contradiction. 
\end{proofclaim}

For each $i \in [k]$, let $U_i$ be the set of vertices which is of colour~$i$ only. 
Suppose that there exists $i \in [k]$ such that $U_i$ is finite, so all but finitely many vertices is of some colour in $[k] \setminus \{i\}$. 
So there exists $j \in [k] \setminus \{i\}$ and a set $A \subseteq \mathbb{N}$ such that each vertex in $A$ is of colour~$j$ and $A$ has upper density at least~$1/(k-1)$. 
By~Claim~\ref{clm:1}, we obtain a monochromatic path~$P$ of colour~$j$ containing~$A$. 
Thus $P$ has upper density at least~$1/(k-1) \ge \phi_{k-1}$ by Proposition~\ref{prop:bipartite-k-colourings-construction}, a contradiction. 
Therefore, $U_i$ is infinite for all $i \in [k]$.

Let $\eps_j = 2^{-j}$ for $j \in \mathbb{N} \cup\{0\}$.
We will now construct monochromatic paths as follows. 

\begin{Claim}
For $j \in \mathbb{N} \cup\{0\}$, there exist an integer~$\ell_j$ and paths $P^j_1,\dots, P^j_k$ such that 
\begin{enumerate}[label = {\rm (\roman*)}]
	\item each $P^j_i$ is monochromatic of colour~$i$ with vertices in~$[\ell_j]$ and both endvertices in~$U_i$; \label{itm:1}
	\item each $P^{j}_i$ is an extension of~$P^{j-1}_i$ and, moreover, $P^{j-1}_i$ is the subgraph of~$P^{j}_i$ induced by the vertex set $[\ell_{j-1}]$; 
	\item there exists $i_j \in [k]$ such that $\frac{| V ( P^j_{i_j} ) \cap [\ell_j] |}{\ell_j} \ge \frac{2(1- \eps_j)}{k+1}$;
	\item $\ell_j \ge \eps_j^{-1} = 2^j$. 
\end{enumerate}
\end{Claim}

\begin{proofclaim}
Let $j=0$. 
Since each $U_i$ is infinite, let $\ell_0$ be the smallest~$\ell$ such that $U_i \cap [\ell] \ne \emptyset$ for all $i \in [\ell]$. 
For each $ i \in [k]$, set $P^0_{i}$ be a vertex in $U_i \cap [\ell_0]$.
Thus the claim holds for $j = 0$. 
Consider $j \in \mathbb{N}$. 
Suppose we have already constructed $P^{j-1}_1,\dots, P^{j-1}_k$ and we construct $P^j_1,\dots, P^j_k$ as follows.

Let $t_0 = t_0(k, \eps_j)$ and $n_0 = n_0(k, \eps_j )$ be given by Corollary~\ref{cor:regularity-ramsey}.
Let $\ell'_j > \ell_{j-1} $ be the smallest integer such that $| U_i \cap [\ell_{j-1}+1, \ell_j'] | \ge t_0$ for all~$i \in [k]$.
Let $m > \ell'_j$ be the smallest integer such that, for all $i \in [k]$, all edges between 
$U_i \cap [\ell'_j]$ and $[m, \infty)$ have colour~$i$.

Let $n' = \max\{ n_0 , 2 m / \eps_j \}$ and $\ell_j = m+ n'$.
Let $I' = [ m+1, \ell_j]$. 
Consider the complete subgraph~$K_{I'}$ of~$K_{\mathbb{N}}$ induced by~$I'$. 
By Corollary~\ref{cor:regularity-ramsey}, $K_{I'}$ contains vertex-disjoint paths~$P_1, \dots, P_t$ such that $t \le t_0$, $\bigcup_{j \in [t]}P_j$ is monochromatic of colour~$i_0$ say, and 
\begin{align*}
\left|\bigcup_{j \in [t]}P_j \right| \ge  \frac{(2-\eps_j) n' }{k+1} \ge \frac{2(1-\eps_j)\ell_j}{k+1}.
\end{align*}
Recall that $ | U_{i_0} \cap [\ell_{j-1}+1, \ell'_j] | \ge t_0 \ge t$ and all edges between $U_{i_0} \cap [\ell'_j]$ and $V(\bigcup_{j \in [t]}P_j )\subseteq [m+1, \infty)$ have colour~$i_0$.
Together with~\ref{itm:1} and using vertices in~$U_{i_0} \cap [\ell_{j-1}+1, \ell'_j]$, we join $P^{j-1}_{i_0},P_1,\dots,P_{t}$ into a monochromatic path~$P^j_{i_0}$ of colour~$i_0$ with endvertices in~$U_{i_0} \cap[\ell'_j]$. 
We are done by setting $P^j_i =P^{j-1}_i$ for $i \in [k] \setminus \{i_0\}$. 
\end{proofclaim}

Note that there exists a colour $i$ such that $i_j = i$ for infinitely many $j \in \mathbb{N}$. 
Then the  monochromatic path $P' = \bigcup_{j \in \mathbb{N}} P^j_i$ satisfies $\overline{d}(P') \ge 2/(k+1) \ge  \phi_{k-1}$, where the last inequality holds by Proposition~\ref{prop:bipartite-k-colourings-construction}.
This is a contradiction.
\end{proof}



\section{Complete bipartite infinite graphs}\label{section:bipartite-work}

We now prove Proposition~\ref{prop:bipartite-k-colourings-construction}, that is, bounding the upper density of monochromatic paths in $k$-edge-coloured $K_{V,W} \in K_{\mathbb{N},\mathbb{N}}$ from above.

\begin{proof}[Proof of Proposition~\ref{prop:bipartite-k-colourings-construction}]
Let $c$ be a proper $k$-edge-colouring of $K_{k,k}$ with vertex classes $X = \{x_1, \ldots, x_k\}$ and $Y= \{y_1, \ldots, y_k\}$.
Let $\phi_X : V \rightarrow X$ and $\phi_Y : W \rightarrow Y$ be such that $\overline{d} (\phi_X^{-1}(x_i) \cup \phi_Y^{-1}(y_j)) = 1/k$ for all $i,j \in [k]$.
(For instance, if $V = \{v_1,v_2, \ldots\}$ and $W= \{w_1,w_2, \ldots\}$ with $v_i < v_{i+1}$ and $w_j < w_{j+1}$, then set $\phi_X(v_i) = x_{i'}$ and $\phi_Y(w_j) = y_{j'}$ such that $i \equiv i' \pmod{k}$ and $j \equiv j' \pmod{k}$.)
We now edge-colour~$K_{V,W}$ such that the edge~$vw$ with $v \in V$ and $w \in W$ has colour~$c(\phi_X(v)\phi_Y(w))$.
Since each colour class of~$c$ is a perfect matching, any monochromatic path in~$K_{V,W}$
lies in $\phi_X^{-1}(x_i) \cup \phi_Y^{-1}(y_j)$ for some $i, j \in [k]$. 
Hence the result follows. 
%
%
\end{proof}

In order to prove Theorem~\ref{theorem:multicolour-bipartite}, we use the notion of an ultrafilter. 
Given an infinite set~$X$, a family~$\mathcal{U}$ of subsets of~$X$ is an \textit{ultrafilter} if $\mathcal{U}$ is closed under finite intersections and supersets, the empty set is not in $\mathcal{U}$, and for every set $Y \subseteq X$, we have that either $Y \in \mathcal{U}$ or $X \setminus Y \in \mathcal{U}$.
Thus if $\mathcal{U}$ is an ultrafilter on~$X$, and $\{X_{1},\ldots, X_{n}\}$ is a finite partition of~$X$, then exactly one $X_{i}$ is in $\mathcal{U}$.
Finally, an ultrafilter~$\mathcal{U}$ is \textit{nonprincipal} if no set in $\mathcal{U}$ is finite.
By Zorn's Lemma, for any infinite set $X$, there exists a nonprincipal ultrafilter on~$X$.



Let $K_{V,W} \in K_{\mathbb{N},\mathbb{N}}$. 
Let $A \subseteq V$ and $B \subseteq W$ be infinite sets. 
A pair $(\mathcal{V}_{A},\mathcal{W}_{B})$ is an \emph{ultrafilter-pair} of~$(A,B)$, if $\mathcal{V}_{A}$ and $\mathcal{W}_{B}$ are nonprincipal ultrafilters on $A$ and~$B$, respectively. 
Given an ultrafilter-pair $\mathcal{S} = (\mathcal{V}_{A},\mathcal{W}_{B})$ of $(A,B)$, define the $k$-vertex-colouring $c_{\mathcal{S}}$ of $V \cup W$ such that $c_{\mathcal{S}}(v) = i$ if $N_i(v) \cap (A \cup B) \in \mathcal{V}_{A} \cup \mathcal{W}_{B}$ for $v \in V \cup W$.
Note that every vertex gets exactly one colour. 
Moreover, for each $i \in [k]$, let
\begin{align}
V_{i}(\mathcal{S}) &= \{v \in V \colon c_{\mathcal{S}}(v) = i\} \text{ and }
W_{i}(\mathcal{S}) = \{w \in W \colon c_{\mathcal{S}}(w) = i\}. \nonumber
\end{align}
When it is clear which ultrafilter-pair we are referring to, we will often omit the $\mathcal{S}$ and write $V_{i}$ and $W_{i}$ instead.  
We make use of the following lemma.

\begin{Lemma}\label{lemma:assisting-lemma}
Let $K_{V,W} \in K_{\mathbb{N},\mathbb{N}}$ be $k$-edge-coloured and $i_0 \in [k]$.
Let $\mathcal{S} = (\mathcal{V}',\mathcal{W}')$ be an ultrafilter-pair on $(V',W')$ with $V' \subseteq V$ and $W' \subseteq W$.
Then there exists a monochromatic path of colour~${i_0}$ containing~$V_{i_0}$.

Moreover, let $U_{i_0}^*$ be the set of vertices $v \in V_{{i_0}}\cup W_{i_0}$ such that $N_{{i_0}}(v) \cap (V_{{i_0}}\cup W_{i_0})$ is infinite. 
If $U_{i_0}^*$ itself is infinite, then there exists a monochromatic path~$P$ of colour~${i_0}$ containing $V_{i_0} \cup W_{i_0}$. 
\end{Lemma}

\begin{proof}
Without loss of generality, we may assume that ${i_0}= 1$.
Let $V_1 = \{a_1, a_2, \dots\}$ with $a_j < a_{j+1}$.
For each $j \in \mathbb{N}$, note that $N_1(a_j) \cap W', N_1(a_{j+1})\cap W' \in \mathcal{W}'$ and recall that $\mathcal{W}'$ is closed under finite intersections, so $N_1(a_j) \cap N_1(a_{j+1})  \cap W' \in \mathcal{W}'$ is infinite. 
Hence, we can find distinct vertices $w_1, w_2, \dots \in W'$ such that $w_j \in N_1(a_j) \cap N_1(a_{j+1})$. 
Then $P = v_1w_1v_2w_2 \dots$ is a monochromatic path of colour~${1}$ containing~$V_{1}$.

We now prove the moreover statement. 
Without loss of generality, we may assume $V^*_1 = V_1 \cap U^*_1$ is infinite. 
Let $V_1 \cup W_1 = \{a_1, a_2, \dots\}$ with $a_j < a_{j+1}$.
Let $A_j = \{a_1, \dots, a_j\}$.
We will construct monochromatic path~$P_{j}$ of colour~$1$ containing $A_j$ with endpoints in~$V^*_1$.
Set $P_0$ be a single vertex in~$V^*_1$.
Suppose that we have constructed~$P_{j-1}$ and construct $P_j$ as follows.
If $a_j \in V(P_{j-1})$, then set $P_j = P_{j-1}$.

Suppose that $a_j \notin V(P_{j-1})$ and $v_{j-1}$ be an endpoint of~$P_{j-1}$.
If $a_j \in V_1$, then pick $v_j \in V^*_1 \setminus (V(P_{j-1}) \cup \{a_j\})$.
Note that $v_{j-1}, v_j, a_j \in V_1$, so $ N_1 (v_{j-1})\cap W'$, $N_1(v_j)\cap W'$ and $N_1(a_j)  \cap W'$ are members of~$\mathcal{W}'$.
Hence, $	N_1 (v_{j-1}) \cap N_1(v_j) \cap N_1(a_j)  \cap W'  \in \mathcal{W}'$ is an infinite set. 
Pick distinct vertices  $w_{j-1}, w_j \in ( N_1 (v_{j-1}) \cap N_1(v_j) \cap N_1(a_j)  \cap W' ) \setminus V(P_{j-1})$.
Note that $P_j = P_{j-1} v_{j-1} w_{j-1} a_j w_j v_j$ is a path of colour~$1$ as desired. 

If $a_j \in W_1$, then pick distinct vertices $v_j \in  V^*_1 \setminus V(P_{j-1})$, $w_{j-1} \in (N_1 (v_{j-1}) \cap W_1)  \setminus V(P_{j-1}) $ and $w_{j} \in (N_1 (v_{j}) \cap W_1)  \setminus V(P_{j-1}) $.
Similarly, pick distinct vertices  $v'_{j-1}, v'_j \in ( N_1 (w_{j-1}) \cap N_1(w_j) \cap N_1(a_j)  \cap V' ) \setminus V(P_{j-1})$.
Note that $P_j = P_{j-1} v_{j-1} w_{j-1} v'_{j-1} a_j v'_j w_j v_j$ is a path of colour~$1$ as desired.
We are done by setting $P = \bigcup_{i \in \mathbb{N}} P_{i}$.
\end{proof}

First we prove Theorem~\ref{theorem:multicolour-bipartite} when $k = 2$. 

\begin{proof}[Proof of Theorem~\ref{theorem:multicolour-bipartite} when $k = 2$]
Fix a $2$-edge-colouring of $K_{V,W}$, and let $\mathcal{S} = (\mathcal{V},\mathcal{W})$ be an ultrafilter-pair on~$(V,W)$.  
Note that $\overline{d}(V_{1} \cup W_{1}) + \overline{d}(V_{2} \cup W_{2}) \ge 1$.
Thus, by relabelling colours if necessary, we may assume that 
\begin{align*}
\overline{d}(V_{1} \cup W_{1}) \ge 1/2.
\end{align*}
We may assume that $\overline{d}(V_{1}), \overline{d}(W_{1}) > 0$ (or else, $\overline{d}(V_{1})\ge 1/2$ or $ \overline{d}(W_{1})\ge 1/2$ and we are done by~Lemma~\ref{lemma:assisting-lemma}).
Hence $V_1$ and $W_1$ are infinite. 

Let $U_{1}^*$ be the set of vertices $v \in V_{1}\cup W_1$ such that $N_{1}(v) \cap (V_1 \cup W_1)$ is infinite.  
If $U^*_1$ is infinite, then Lemma~\ref{lemma:assisting-lemma} implies that there is a path of colour~$1$ containing $V_1 \cup W_1$, as required. 
Thus we may assume that $U^*_1$ is finite.

Thus $V_1 \setminus U^*_1$ and $W_1 \setminus U^*_1$ are infinite.
Futhermore, every $v \in V_1 \setminus U^*_1$ (and $w \in W_1 \setminus U^*_1$) sends finitely many edges of colour~$1$ to~$W_1$ (and to~$V_1$, respectively). 
It is easy to construct a monochromatic path~$P$ of colour~$2$ with vertex set $V(P) = (V_1 \cup W_1 ) \setminus U^*_1$ (see the proof of the moreover statement of Lemma~\ref{lemma:assisting-lemma}).
Note that $\overline{d}(P) = \overline{d}((V_1 \cup W_1 ) \setminus U^*_1) = \overline{d}(V_{1} \cup W_{1}) \ge 1/2$ as required.
%
\end{proof}

Before proving Theorem~\ref{theorem:multicolour-bipartite} when $k \ge 3$, we would need to define the lower density of a set. 
Given a set $A \subseteq \mathbb{N}$, the \textit{lower density of~$A$} is defined as 
\begin{equation}
\underline{d}(A) = \liminf_{n \rightarrow \infty}\frac{|A \cap [n]|}{n}. \nonumber
\end{equation}
For sets $U, W \subseteq \mathbb{N}$ with $U \cap W$ finite (i.e. almost disjoint sets of $U, W \subseteq \mathbb{N}$), the following standard inequality holds:
\begin{align}
\underline{d}(U)+ \underline{d}(W)
\le
\underline{d}(U \cup W)
\le 
\underline{d}(U) + \overline{d}(W)
\le 
\overline{d}(U \cup W)
\le 
\overline{d}(U) + \overline{d}(W)
.
\label{eqn:density}
\end{align}

\begin{proof}[Proof of Theorem~\ref{theorem:multicolour-bipartite} when $k \ge 3$]
Let $\alpha = 1/(2k-3)$.
Suppose  to the contrary that there exists $\eps>0$ and a $k$-edge-coloured $K_{V,W}$ such that every monochromatic path~$P$ in $K_{V,W}$ has $\overline{d}(P) \leqslant \alpha - 2\eps$.

Let $\beta$ be the supremum of $\max_{i \in [k]}\{\overline{d}(V_{i}(\mathcal{S}) \cup W_{i}(\mathcal{S})) \}$ taken over all ultrafilter-pairs~$\mathcal{S}$ of~$(V',W')$ with infinite sets $V' \subseteq V$ and $W' \subseteq W$.
Clearly, 
\begin{align}
	\beta \ge 1/k. \label{eqn:beta}
\end{align}
Let $\mathcal{S}^{0} = ( \mathcal{V}^{0},\mathcal{W}^{0})$ be an ultrafilter-pair on $(V^0,W^0)$ with infinite sets $V^0 \subseteq V$ and $W^0 \subseteq W$ such that $\max_{i \in [k]}\{\overline{d}(V_{i}(\mathcal{S}^{0}) \cup W_{i}(\mathcal{S}^{0})) \} \ge \beta - \eps$. 
For each $i \in [k]$, let $V_{i}^{0} = V_{i}(\mathcal{S}^{0})$ and $W^0_i = W_{i}(\mathcal{S}^{0})$.
By relabelling if necessary, we may assume that 
\begin{align}
\label{equation:1}
\overline{d}(V_{1}^{0} \cup W_{1}^{0})  = \max_{i \in [k]}\{\overline{d}(V_{i}^{0} \cup W_{i}^{0}) \} \ge \beta - \eps. 
\end{align}
By Lemma~\ref{lemma:assisting-lemma} (with $\mathcal{S} = \mathcal{S}^0$ and $i_0 =1$), there exists a path of colour~$1$ containing $V_{1}^0$, so $\overline{d}(V_{1}^{0}) \le \alpha - 2\eps$.
Hence
\begin{align*}
\overline{d}(W_1^0) \overset{\text{\eqref{eqn:density}}}{\ge} 
\overline{d}(V_{1}^0 \cup W_1^0) - \overline{d}(V_{1}^{0}) 
\overset{\text{\eqref{equation:1}}}{\ge}
 \beta - \alpha + \eps \overset{\text{\eqref{eqn:beta}}}{>}0
\end{align*}
 and so $W^0_1$ is infinite.
Similarly, $\overline{d}(W_{1}^{0}) \le \alpha - 2 \eps$ and $V^0_1$ is infinite. 
Moreover,
\begin{align}\label{equation:4}
\beta \overset{\mathclap{\text{\eqref{equation:1}}}}{\le} \overline{d}(V_{1}^{0} \cup W_{1}^{0}) + \eps \overset{\mathclap{\text{\eqref{eqn:density}}}}{\le}
 \overline{d}(V_{1}^{0}) + \overline{d}(W_{1}^{0}) + \eps \le 2\alpha - 3\eps < 2 \alpha.
\end{align}

Let $\mathcal{S}'= (\mathcal{V}',\mathcal{W}')$ be an ultrafilter-pair on $(V_{1}^{0},W_{1}^{0})$.  
Let $V'_i = V_i (\mathcal{S}')$ and $W'_i = W_i (\mathcal{S}')$ for all $ i \in [k]$.

\begin{Claim} \label{clm:finite}
$(V'_1 \cup W'_1) \cap (V^0_1 \cup W^0_1) $ is finite. 
\end{Claim}

\begin{proofclaim}
Suppose to the contrary that $(V'_1 \cup W'_1) \cap (V^0_1 \cup W^0_1) $ is infinite.
Without loss of generality $V_{1}^{0} \cap V_{1}'$ is infinite.
For all $v \in V_{1}^{0} \cap V_{1}' \subseteq  V_{1}'$, we have that $N_1(v) \cap W^0_1 \in \mathcal{W}'$ is infinite. 
Lemma~\ref{lemma:assisting-lemma} (with $\mathcal{S} = \mathcal{S}$, $i_0 = 1$ and $U^*_1 \supseteq V_{1}^{0} \cap V_{1}'$) implies that there exists a path of colour~1 containing $V_{1}^{0} \cup W_{1}^{0}$ with upper density at least
\begin{align*}
 \overline{d}(V_{1}^{0} \cup W_{1}^{0})  \ge \beta - \eps \overset{\eqref{eqn:beta}}{\ge} 1/k - \eps \ge \alpha - \eps,
\end{align*}
a contradiction. 
\end{proofclaim}

Consider the ultrafilter-pair $\mathcal{S}^* = ( \mathcal{V}', \mathcal{W}^0)$, so $V_1(\mathcal{S}^*) = V_1^0$ and $W_1(\mathcal{S}^*) = W'_1$. 
Moreover, for all $w \in W_{1}'$, $N_1(w) \cap V_{1}^{0} \in \mathcal{V}'$ is infinite. 
If $W_{1}'$ is finite, then $\overline{d}(V_{1}^{0} \cup W_{1}') = \overline{d}(V_{1}^{0})\le \alpha -2 \eps$.
If $W_{1}'$ is infinite, then Lemma~\ref{lemma:assisting-lemma} (with $\mathcal{S} = \mathcal{S}^*$, $i_0 = 1$ and $U^*_1 \supseteq W_1'$) implies that there exists a path of colour~$1$ containing $V_1^0 \cup W_{1}'$. 
In both cases, we have $\overline{d}(V_{1}^{0} \cup W_{1}') \le \alpha -2 \eps$.
Similarly by considering the ultrafilter-pair $(\mathcal{V}^0,\mathcal{W}')$, we deduce that $\overline{d}(W_{1}^{0} \cup V_{1}') \le \alpha -2 \eps$. 
Together with Claim~\ref{clm:finite} and \eqref{eqn:density}, we deduce that 
\begin{align*}
\underline{d}(V_{1}' \cup  W_{1}') + \beta - \eps
& \overset{\mathclap{\eqref{equation:1}}}{\le}
\underline{d}(V_{1}' \cup W_{1}') + \overline{d}(V_{1}^0 \cup W_{1}^0) 
\\ &
\le
 \overline{d}(V_{1}' \cup  W_{1}' \cup V_{1}^0 \cup W_{1}^0 ) 
\\&
\le \overline{d}(V_{1}^{0} \cup W_{1}') +\overline{d}(W_{1}^{0} \cup V_{1}') 
\le 2(\alpha - 2 \eps),
\end{align*}
which implies that 
\begin{align} 
	\underline{d}(V_{1}' \cup W_{1}') & \le  2\alpha - \beta.
\label{equation:2}
\end{align}

Since $\mathcal{V}'$ is an ultrafilter on~$V_{1}^{0}$, there exists some $i_0 \in [k]$ such that $V_{1}^{0} \cap V_{i_0}' \in \mathcal{V}'$.
Note that $V_{1}^{0} \cap V_{i_0}'$ is infinite, so $i_0 \ne 1$ by Claim~\ref{clm:finite}.
Therefore, without loss of generality, we may assume that $i_0 =2$. 
Recall that $V_{1}^{0} \cap V_{2}' \in \mathcal{V}'$.
For all $w \in W_2'$, note that $N_2(w) \cap V^0_1 \in \mathcal{V}'$ implying that $N_2(w) \cap V^0_1 \cap V_{2}' = (N_2(w) \cap V^0_1) \cap (V_{1}^{0} \cap V_{2}')\in \mathcal{V}'$ is infinite. 
If $W_2'$ is infinite, then Lemma~\ref{lemma:assisting-lemma} (with $\mathcal{S} = \mathcal{S}'$, $i_0 = 2$ and $U^*_2 \supseteq W_2'$) implies that there exists a path of colour~$2$ containing $V_{2}' \cup W_2'$.
If $W_2'$ is finite, then we have $\overline{d}(V_{2}' \cup W_{2}') = \overline{d}(V_{2}' )$ and Lemma~\ref{lemma:assisting-lemma} implies that  there is a path of colour~2 containing~$V_{2}'$.
In both cases, we deduce that
\begin{align} \label{equation:3}
\overline{d}(V_{2}' \cup W_{2}') \le \alpha - 2\eps < \alpha . 
\end{align}
Recall the definition of~$\beta$ that $\overline{d}(V_{i}' \cup W_{i}') \le \beta$ for all $i \in [k]$. 
Putting these all together, we get that
\begin{align*}
1 & = \underline{d}(V \cup W) \le \underline{d}(V_{1}' \cup W_{1}') + \overline{d} \left(\bigcup_{2 \le i \le k}(V_{i}' \cup W_{i}') \right)  
  \le \underline{d}(V_{1}' \cup W_{1}') + \sum_{ 2 \le i \le k } \overline{d}(V_{i}' \cup W_{i}')  \\
	&
\overset{\mathclap{\text{\eqref{equation:2}, \eqref{equation:3}}}}{<}
  ( 2 \alpha - \beta) + \alpha + \sum_{3 \le i \le k} \overline{d}(V_{i}' \cup W_{i}') 
\le  3 \alpha - \beta + (k-2)\beta  
= 3\alpha +(k-3)\beta  \\
&\overset{\mathclap{\text{\eqref{equation:4}}}}{<} (2k-3)\alpha = 1,
\end{align*}
a contradiction.
\end{proof}

\section*{Acknowledgments}
The authors would like to thank Louis DeBiasio and the referees for their valuable comments
and for their constructive suggestions on the presentation of the paper.
In particular, we are grateful for one referee for improving our result of Theorem~\ref{theorem:multicolour-bipartite}.

\end{document}